\documentclass[preprint]{elsarticle}

\usepackage{amsmath}
\usepackage{amsthm}
\usepackage{amsfonts}
\usepackage{amssymb}
\usepackage{bm}
\usepackage{mathrsfs}
\usepackage{mathtools}
\usepackage{curves}
\usepackage{subfig}
\usepackage{booktabs}
\usepackage{tabularx}
\usepackage{multirow}
\usepackage{algorithm}
\usepackage{algorithmic}
\usepackage{listings}
\usepackage{tikz}
\usepackage{tikz-3dplot}
\usepackage{color}
\usepackage{animate}
\usepackage{epstopdf}
\usepackage{bookmark}
\usepackage{caption}
\usepackage{enumitem}
\usepackage[normalem]{ulem}
\usepackage[capitalise, nameinlink]{cleveref}

\usetikzlibrary{shapes.geometric, arrows, fit}

\usepackage{tensor}      
\newcommand{\fd}[4][]{\tensor*[^{#1}_{#3}]{\mathcal{D}}{_{#4}^{#2}}}
\newcommand{\fri}[4][]{\tensor*[^{#1}_{#3}]{I}{_{#4}^{#2}}}

\newcommand{\rfd}[3]{\fd[\text{\tiny RL}]{#1}{#2}{#3}}

 \usepackage{diagbox}
\usepackage{geometry}
\usepackage{natbib}
\biboptions{sort&compress}
\graphicspath{{figures/}}

\theoremstyle{definition}

\newdefinition{definition}{Definition}
\newdefinition{example}{Example}[section]

\journal{}

\begin{document}
\begin{frontmatter}
	\title{Finite Element Methods for Fractional PDEs in Three Dimensions}
	\author{Zongze Yang}
	\author{Yufeng Nie}
	\author{Zhanbin Yuan}
	\author{Jungang Wang\corref{cor1}}\ead{wangjungang@nwpu.edu.cn}
			
	\cortext[cor1]{Corresponding author}
			
	\address{Department of Applied Mathematics,
		Northwestern Polytechnical University,
	Xi'an, Shaanxi 710072, China}
			
	\begin{abstract}
		This paper is a generalization of the previous work (Yang et.al, J. Comput. Phys. 330 (2017), 863$-$883) to the 3-D irregular convex domains. The analytical calculation formula of fractional derivatives of finite element basis functions are given and a path searching method is developed to find the integration paths corresponding to the Gaussian points. Moreover, a template matrix is introduced to speed up the procedures. Numerical experiments on a steady problem are presented verifying the efficiency of proposed techniques.
	\end{abstract}
			
	\begin{keyword}
		finite element method \sep
		fractional derivatives \sep
		three dimensions \sep
		irregular domains  \sep
		path searching method \sep
	\end{keyword}
\end{frontmatter}

\section{Introduction}\label{sec:introduction}
Fractional calculus has been widely used as a modeling tool in numerous fields for diverse applications. However, due to their non-local properties, it is difficult to get closed form or numerical solutions of fractional partial differential equations, especially on the irregular domains in high dimensions. Based on the theoretical and computational framework \cite{ervinVariationalFormulationStationary2006,ervinVariationalSolutionFractional2007,roopComputationalAspectsFEM2006} developed by Ervin and Roop, {Zhao et al.\cite{zhao2015} developed finite element method for two-dimensional space-fractional advection-dispersion equations}. Yang et al.\cite{yangFiniteElementMethod2017a}  proposed a finite element method(FEM) on unstructured meshes  for nonlinear Riesz fractional partial differential equations(FPDEs) in {two-dimensional(2-D)} domains. Furthermore, Fan et al.\cite{fanNovelUnstructuredMesh2017a,fanEfficientFiniteElement2018} solved the time-space fractional wave equation and fractional Schr\"odinger equation using unstructured meshes. Dehghan and Abbaszadeh \cite{dehghanEfficientTechniqueBased2018,dehghanFiniteDifferenceFinite2018} studied the finite element solutions for space-multi-time fractional Bloch-Torrey equations and tempered fractional diffusion-wave equations. Though there are so many works on FEM for FPDEs, most of them focus on 1-D or 2-D problems. FPDEs on {three-dimensional(3-D)} domains are seldom solved by FEM in literature. The main purpose of this work is to extend FEM on unstructured meshes to solve FPDEs on 3-D irregular domains.

Let $\mathbf{x} = (x_1, \cdots, x_n)^T$, $\mathbf{y} = (y_1, \cdots, y_n)^T$ be points in $\mathbb{R}^n$ and $\mathbf{e}_i (i=1, 2, \cdots, n)$ be the $i$-th unit column vector {where $n = 2$ or $3$}. The boundary of a convex domain $\varOmega$ can be characterized by $a_i(\mathbf{x})$ and $b_i(\mathbf{x})$ as below,
\begin{equation}
	\begin{aligned}
		a_i(\mathbf{x}) & := \inf\{y_i| \mathbf{y} = \mathbf{x} + k \mathbf{e}_i, 
		\mathbf{y}\in \varOmega, k \in \mathbb{R} \},                           \\
		b_i(\mathbf{x}) & := \sup\{y_i| \mathbf{y} = \mathbf{x} + k \mathbf{e}_i, 
		\mathbf{y}\in \varOmega, k \in \mathbb{R} \}.                           
	\end{aligned}
\end{equation}
Here, $a_i(\mathbf{x})$ and $b_{i}(\mathbf{x})$ respectively represent the lower and upper bounds of
the line segment parallel to $\mathbf{e}_i$ through point $\mathbf{x}$.

\begin{definition}
	Let $\alpha > 0$,  $f(\mathbf{x}) \in L_1(\varOmega)$. Then the {left and 
  right} Riemann-Liouville fractional
	integral operators $\fri{\alpha}{a_i(\mathbf{x})}{x_i}$ and $\fri{\alpha}{x_i}{b_i(\mathbf{x})}$
  of order $\alpha>0$ {with respect to $x_i$} are defined by
	\begin{equation}\label{fcl}
		\fri{\alpha}{a_i(\mathbf{x})}{x_i} f(\mathbf{x})
		:= \frac{1}{\Gamma(\alpha)} \int_{a_i(\mathbf{x})}^{x_i}(x_i-y_i)^{\alpha-1}
		f\big(\mathbf{x} + (y_i - x_i)\mathbf{e}_i\big)d y_i,
	\end{equation}
	\begin{equation}\label{fcr}
		\fri{\alpha}{x_i}{b_i(\mathbf{x})} f(\mathbf{x})
		:= \frac{1}{\Gamma(\alpha)}\int_{x_i}^{b_i(\mathbf{x})}(y_i-x_i)^{\alpha-1}
		f\big(\mathbf{x} + (y_i - x_i)\mathbf{e}_i\big)d y_i.
	\end{equation}
\end{definition}

\begin{definition}
  Let $m-1 <\alpha < m, m\in \mathbb{N}_+$. Then the {left and right}
  Riemann-Liouville fractional derivatives of {$f(\mathbf{x})$}	of order $\alpha$ 
  {with respect to $x_i$}
  are defined by
	\begin{equation}\label{def}
		\begin{aligned}
			\rfd{\alpha}{a_i(\mathbf{x})}{x_i} f(\mathbf{x}) & :=\frac{d^{m}}{dx_i^{m}}       
			\fri{m-\alpha}{a_i(\mathbf{x})}{x_i} f(\mathbf{x}), \\
			\rfd{\alpha}{x_i}{b_i(\mathbf{x})} f(\mathbf{x}) & :=(-1)^m\frac{d^{m}}{dx_i^{m}} 
			\fri{m-\alpha}{x_i}{b_i(\mathbf{x})} f(\mathbf{x}).
		\end{aligned}
	\end{equation}
\end{definition}

Fractional differential operators in application mainly contain three forms: 
{
$ c(\mathbf{x}) \mathcal{D}^{2\mu}_i u(\mathbf{x}) $, 
$\mathcal{D}_i \big(c(\mathbf{x}) \mathcal{D}^{\nu}_i u(\mathbf{x}) \big)$ 
or $\mathcal{D}^{\nu}_i \big( c(\mathbf{x}) \mathcal{D}_i u(\mathbf{x}) \big)$, 
where $1/2 < \mu < 1$, $0 < \nu < 1$, and $\mathcal{D}^{\gamma}_i (\gamma > 0)$ 
}
is left or right
Riemann-Liouville fractional derivative operator with respect to $x_{i}$, and $\mathcal{D}_{i}$ is
an {abbreviation} of classical first order derivative operator $\frac{\partial}{\partial x_i}$.
When $u(\mathbf{x})$ is zero on the boundary, the variational forms of those three cases have a uniform form 
{
\begin{equation}\label{equ:varsss}
	\big( c(\mathbf{x}) \mathcal{D}^{\alpha}_i u(\mathbf{x}),\mathcal{\hat D}^{\beta}_i v(\mathbf{x}) \big),
\end{equation}
where $ \alpha,\beta \ge 0$, $\alpha + \beta = 2\mu$ (or $\alpha + \beta = 1 + \nu$ for the last two forms), and $\mathcal{\hat D}^{\beta}_i$ represents
the right (left) fractional operator with respect to $x_i$ if
$\mathcal{D}^{\alpha}_i$ represents the left (right) fractional operator.
}

In this paper, we mainly discuss how to assemble stiffness matrix for \eqref{equ:varsss}. From \eqref{def}, we know that the calculation of fractional derivatives of finite element basis functions contain integrals over the lines connecting the Gaussian points to the boundary along the $i-$th direction which we call the {\it {integration path}} in the following. This is the main challenge in the implementation of FEM for FPDEs and will be discussed in detail in next section. The second section also includes finding the {integration paths} and some speed-up {techniques}. In the last section, the proposed methods are applied to solve 3-D steady fractional problems.

\section{{Assembling the fractional stiffness matrix of $\big( c(\mathbf{x}) \mathcal{D}^{\alpha}_i u(\mathbf{x}),\mathcal{\hat D}^{\beta}_i v(\mathbf{x}) \big)$ }}
Assembling the fractional stiffness matrix is the key point to solve FPDEs using FEM.
In this section, we first deduce the analytical calculation formula of fractional derivatives of finite element basis functions. Then, we present the method of finding the integration paths for Gaussian points. Finally, some speed-up {techniques} are discussed.

\subsection{{Deduction of the analytical formula of fractional derivative of FE basis functions}}\label{ss:bp}
Assume the domain $\varOmega \in \mathbb{R}^{n}$ is a polygonal domain which can be partitioned into
simplexes in $\mathbb{R}^{n}$.  Let ${\{\mathcal{T}_h\}}$ be a family of  regular
partitions of $\varOmega$, and $h$ be the maximum diameter of elements in $\mathcal{T}_h$.  For
finite element methods, the idea is to approximate solutions of equations in  finite dimensional spaces.  So we define the test and trial space $V_h = \{v_h: v_h \in C(\varOmega),v_h \in V, v_h|_E\in P_s(E), \forall E \in \mathcal{T}_h\}$,
where $P_s(E)$ is the set of polynomials of degree $\leq s$ in $E$. Denote the basis functions of $V_h$ as $\psi_l(\mathbf{x}), l = 1, 2, \cdots, \mathcal{N}$. The stiffness matrix can be calculated by adding the element stiffness matrices, i.e.,
\begin{equation}
	K = \sum_{E\in\mathcal{T}_h} K_E,
\end{equation}
where $K_E$ is element stiffness matrix with $(K_E)_{kl} = \big( c(\mathbf{x}) \mathcal{D}^{\alpha}_i \psi_l(\mathbf{x}),{\mathcal{\hat D}^{\beta}_i}\psi_k(\mathbf{x}) \big)_E$ in which $(\cdot,\cdot)_E$ represents the $L^2$ inner product on element $E$.

By Gaussian quadrature on $E$, we can approximate $(K_E)_{kl}$ by
\begin{equation}
  (K_E)_{kl} = \sum_{j=1}^m w_j c(\mathbf{x}_j) \mathcal{D}^{\alpha}_i \psi_l(\mathbf{x}_j) {\mathcal{\hat D}^{\beta}_i} \psi_k(\mathbf{x}_j),
\end{equation}
where $\mathbf{x}_j$ and $w_j$ $(j=1,2,\cdots, m)$ are points and weights of Gaussian quadrature rule on $E$.

\begin{figure}[h]
	\centering
	\subfloat[]{\label{path2d}\includegraphics[width=.45\textwidth]{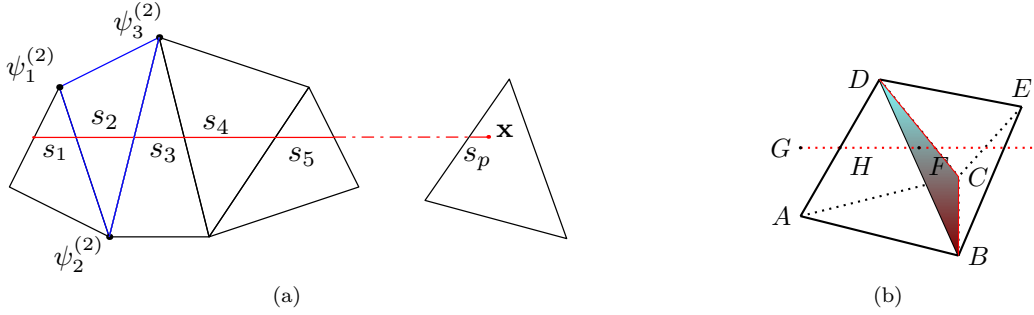}}
	\qquad\qquad\qquad
	\subfloat[]{\label{path3d}
		\tdplotsetmaincoords{120}{90}
		\centering
		\begin{tikzpicture}[scale=2.1, tdplot_main_coords]
			          
			\coordinate (A) at (0, 0, 0);
			\coordinate (B) at (-0.5, 1, 0);
			\coordinate (C) at (0.5, 1, 0);
			\coordinate (D) at (0, 0.5, 1);
			\coordinate (E) at (0, 1.4, 0.8);
			\coordinate (G) at (0, 0, 0.5);
			\coordinate (F) at (0, 0.75, 0.5);
			\coordinate (I) at (0, 0.25, 0.5);
			\coordinate (H) at (0, 1.5, 0.5);
			\draw[-, thick]   (A) --  (B)
			(B) --  (E)
			(E) --  (D)
			(B) --  (D)
			(D) --  (A);
			\draw[dotted, thick]  (A) --  (C) (E) --  (C) (B) --  (C) (D) --  (C);
			\filldraw[draw=red,bottom color=red!50!black, top color=cyan!50]  (B) -- (C)  -- (D);
			\draw[dotted,  thick, red] (G) -- (F) -> (H);
			          
			\fill[black]  (F) circle [radius=0.3pt];
			\fill[black]  (G) circle [radius=0.3pt];
			\fill[black]  (I) circle [radius=0.3pt];
			          
			\draw (A) node [left]            {$A$}
			(B) node [right]           {$B$}
			(C) node [right]           {$C$}
			(D) node [left]            {$D$}
			(F) node [below right]     {$F$}
			(G) node [left]            {$G$}
			(E) node [above]           {$E$}
			(I) node [below right]     {$H$};
			          
		\end{tikzpicture}
	}
	\caption{(a). Path of integration for Gaussian point $\mathbf{x}$ on 2-D domain. $s_k (k=1,2, \cdots, p)$ are the segments of the integral path.  $\psi_1^{(2)}, \psi_2^{(2)}, \psi_3^{(2)}$ are local basis
		functions in $E_2$.\qquad (b). Path of integration across two elements in 3-D.}\label{path}
\end{figure}

We take the left Riemann-Liouville fractional operator as an example and suppose we have found the integration path for Gaussian point $\mathbf{x}$. Denote $s_k (k=1,2, \cdots, p)$ the segments of the integral path, see \cref{path}. And assume $s_k$ is
in element $E_k$.  Then we have
\begin{equation}
	\begin{aligned}
		\mathcal{D}^{\alpha}_i \psi_l(\mathbf{x})                              
		= \frac{1}{\Gamma(1-\alpha)} \frac{d}{d x_i}                           
		\int_{a_i(\mathbf{x})}^{x_i}(x_i-y_i)^{-\alpha}\psi_l(\mathbf{y})d y_i 
		= \sum_{k=1}^{p} \mathcal{D}^{\alpha}_{i, s_k} \psi_l(\mathbf{x}),     
	\end{aligned}
\end{equation}
where $\mathbf{y} = \mathbf{x} + (y_i - x_i)\mathbf{e}_i$ and 
\begin{equation}
	\mathcal{D}^{\alpha}_{i, s_k} \psi_l(\mathbf{x})
	:= \frac{1}{\Gamma(1-\alpha)}\frac{d}{d x_i} \int_{s_k}(x_i-y_i)^{-\alpha}\psi_l(\mathbf{y})d y_i.
\end{equation}

From the definition of $\mathcal{D}^{\alpha}_{i, s_k}\psi_l(\mathbf{x})$, we know that
{\it if $\psi_l$ is not a local basis function on $E_k$, i.e. $\psi_l$ is zero on $s_k$,
	we have $\mathcal{D}^{\alpha}_{i, s_k}\psi_l(\mathbf{x}) = 0$.} So, we only need to compute
$\mathcal{D}^{\alpha}_{i, s_k}\psi_l^{(k)}(\mathbf{x})$ where $\psi_l^{(k)} (l = 1, 2, \cdots, n+1)$
are local basis functions on $E_k$.  In this way, we can calculate the element stiffness matrix locally.

\subsection{Searching the integration paths for Gaussian points}\label{ss:fip}
The integration paths consist of intersection segments of elements with the rays starting from Gaussian points along a certain direction. Most existing algorithms for searching the integration paths are exhaustive, i.e.~computing the intersection segments of the ray with most of the elements, even every element in the mesh. In
\cite{yangFiniteElementMethod2017a}, we introduced the effect domain to reduce the computation time. Here, we provide another way to accelerate the computation by the ray-simplex intersection algorithm. The key idea of this algorithm is to represent the ray in a new coordinate system related to the simplex \cite[Chapter~2]{glassner1989introduction}.

Denote the ray and the simplex as $R$ and $S$, separately. Assume the vertexes of simplex $S$ are
$\mathbf{v}_i (i=0,\cdots, n)$ and set
$\tilde{\mathbf{e}}_i = \mathbf v_i - \mathbf v_0\, (i = 1, \cdots, n)$. It is easy to see
that $\{\tilde{\mathbf{e}}_{i}\}_{i=1}^n$ forms a basis for $\mathbb{R}^{n}$. Hence, every point in
simplex $S$ can be expressed as
\begin{equation}
	\mathbf{x} = \mathbf{v}_0 + k_1 \tilde{\mathbf{e}}_1 + k_2 \tilde{\mathbf{e}}_2
	+ \cdots + k_n \tilde{\mathbf{e}}_n,
\end{equation}
where 
\begin{equation}\label{eq:ck}
	\left\{
	\begin{aligned}
		  & k_i \ge 0,                                  \\
		  & k_0 := 1 - \sum\nolimits_{i=1}^n k_i \ge 0. 
	\end{aligned}\right.
\end{equation}
If we set $\tilde{\mathbf{k}} = (k_0, k_1, \cdots, k_n)^T$, then $\tilde{\mathbf{k}}$ is the volume
coordinate of $\mathbf x$ in the simplex.  Let $\mathbf{u}_0$ be the start point of $R$ and $\mathbf
d (|\mathbf d| = 1)$ be the direction of $R$, then every point on ray $R$ can be written as
\begin{equation}
	\mathbf{x} =  \mathbf{u}_0 + r \mathbf{d},\quad r \ge 0.
\end{equation}
Let $ A = (\tilde{\mathbf{e}}_1, \tilde{\mathbf{e}}_2, \cdots, \tilde{\mathbf{e}}_n)$,
$\mathbf{ k} = (k_1, k_2, \cdots, k_n)^T$, $\mathbf{b} = \mathbf{u}_0 - \mathbf v_0$.
If $\mathbf x$ is an intersection point of $S$ and $R$, we have
\begin{equation*}\label{eq:condrs}
	A\mathbf{ k} = \mathbf{b} + r \mathbf{d},\quad
	r \ge 0, \quad k_i \ge 0 \quad (i=1, 2, \cdots, n),
\end{equation*}
i.e.
\begin{equation}\label{eq:condrs2}
	\mathbf{ k} = A^{-1}\mathbf{b} + r A^{-1}\mathbf{d},\quad
	r \ge 0, \quad k_i \ge 0 \quad (i=1, 2, \cdots, n).
\end{equation}
Inserting \eqref{eq:condrs2} into \eqref{eq:ck}, we obtain the intersection conditions 
including only unknown variable $r$ as below
\begin{equation}\label{eq:cond}
	\left\{
	\begin{aligned}
		  & r \ge 0,                                                  \\
		  & \tilde{\mathbf{b}} + r \tilde{\mathbf{d}} \ge \mathbf{0}, 
	\end{aligned}\right.
\end{equation}
where the operator `$\ge$' is element-wise and
\begin{equation}
	\tilde{\mathbf{b}} = \left(
	\begin{array}{c}
		A^{-1}\mathbf{b}                   \\
		1-\mathbf{e}^{T}  A^{-1}\mathbf{b} 
	\end{array}
	\right) ,\quad
	\tilde{\mathbf{d}} = \left(
	\begin{array}{c}
		A^{-1}\mathbf{d}                 \\
		-\mathbf{e}^{T} A^{-1}\mathbf{d} 
	\end{array}
	\right),
\end{equation}
where $\mathbf e= (1, 1, \cdots, 1)^T$.

Now, the intersection segment can be derived by solving inequalities \eqref{eq:cond}.  Actually, we
only need to know the maximum and the minimum $r$ satisfying \eqref{eq:cond}.  Let $I$ be the set of
all $r$ satisfying \eqref{eq:cond}. If $I$ is empty or contains only one point, $I$ has no contribution to the integral.
Otherwise, $I$ is an interval, i.e.~$I = [r_{\min}, r_{\max}], r_{\min} < r_{\max}$. Therefore, the
two end points of the intersection segment are
\begin{equation}
  \mathbf x_{0} = {\mathbf u_{0}} + r_{\min}\mathbf d, \quad
  \mathbf x_{1} = {\mathbf u_{0}} + r_{\max}\mathbf d,
\end{equation}
where $\mathbf{x}_{1}$ is the go-out point of $R$ from simplex $S$.
The volume coordinates of $\mathbf{x}_{0}$ and $\mathbf{x}_{1}$ can be written as
\begin{equation}\label{equ:ac}
	\tilde{\mathbf{k}}_{\min} = \tilde{\mathbf{b}} + r_{\min} \tilde{\mathbf{d}}, \quad
	\tilde{\mathbf{k}}_{\max} = \tilde{\mathbf{b}} + r_{\max} \tilde{\mathbf{d}}.
\end{equation}

\begin{figure}[!h]
	\centering
	\begin{minipage}{0.5\linewidth}\noindent
		\begin{tikzpicture}[scale=3]
			\coordinate (A) at (0.5, 0.75, 0.433);
			\coordinate (B) at (0, 0, 0);
			\coordinate (C) at (0.5, 0, 0.866);
			\coordinate (D) at (1, 0, 0);
									
			\node [above] at (A) {$A$ $(1, 0, 0, 0)$};
			\node [left]  at (B) {$B$ $(0, 1, 0, 0)$};
			\node [below] at (C) {$C$ $(0, 0, 1, 0)$};
			\node [right] at (D) {$D$};
			\node [below, xshift=0.5cm, yshift=-0.2cm] at (D) {$(0, 0, 0, 1)$};
									
			\draw (A) -- (B) -- (C) -- (D) -- (A);
			\draw (A) -- (C);
			\draw[dotted] (B) -- (D);
			\filldraw[draw=red,bottom color=red!50!black, top color=cyan!50] (A) -- (B)  -- (C);
			\filldraw[draw=red,bottom color=red!50!black, top color=cyan!50] (A) -- (C)  -- (D);
		\end{tikzpicture}
	\end{minipage}
	\begin{minipage}{0.48\linewidth}
		\begin{tabular}{cc}
			\toprule
			Position of    & Volume                \\ 
			Points         & Coordinates           \\
			\midrule
      $AB$           & {$(k_0, k_1, 0, 0)$}.   \\
      $AC$           & {$(k_0, 0, k_2, 0)$}.   \\
      $AD$           & {$(k_0, 0, 0, k_3)$}.   \\
      $BC$           & {$(0, k_1, k_2, 0)$}.   \\
      $CD$           & {$(0, 0, k_2, k_3)$}.   \\
      $\Delta_{ABC}$ & {$(k_0, k_1, k_2, 0)$}. \\
      $\Delta_{ABD}$ & {$(k_0, k_1, 0, k_3)$}. \\
      $\Delta_{ACD}$ & {$(k_0, 0, k_2, k_3)$}. \\
      $\Delta_{BCD}$ & {$(0, k_1, k_2, k_3)$}. \\
			\bottomrule
		\end{tabular}
	\end{minipage}
	\caption{The volume coordinates of nodes and points on edges and faces in tetrahedron.}
	\label{fig:vc}
\end{figure}
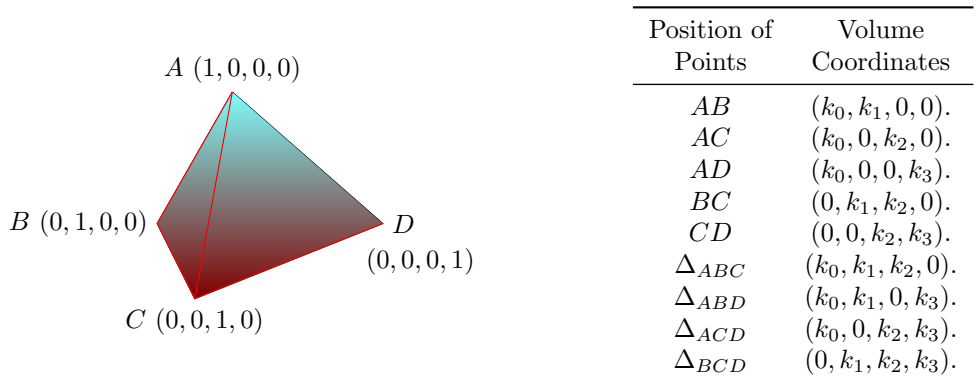

Assume we have known the adjacent information of the triangulation. If a ray intersects with a
simplex $\Delta$, then it must go out from some point on some $m$-face
($0 \le m < n, m \in \mathbb{N}$) of the simplex $\Delta$, for example, edge(1-face), vertex(0-face)
in 2-D case, and face(2-face), edge(1-face), vertex(0-face) in 3-D case.  By the volume coordinate of the
go-out point, we can get the simplex with which the ray intersect after going out of
current simplex.
{\it If the volume coordinate of one point has $k$ zero components, then the point must be on some
	$(n-k)$-face of the simplex. And the exact face can be determined by the distribution of the zero
	components of the volume coordinate, see \cref{fig:vc} for 3-D case.}
If the go-out point is on an $m$-face, we can search all the elements (simplex) which contain the
$m$-face. Usually, most of the go-out points will locate on $(n-1)$-face of the simplex $\Delta$,
for example, edges in 2-D and faces in 3-D (See \cref{path3d}). And in this situation, there are at
most one simplex shared the face with simplex $\Delta$. This is why this algorithm is of high-efficiency.

\subsection{Speedup addition of sparse matrix in MATLAB}\label{sss:imp}
When assembling the stiffness matrix we need to add element stiffness matrices together.  However,
this operation will cost much CPU time because of the sparse matrix storage method in MATLAB.
In this subsection, we will exploit the sparsity of the fractional stiffness matrix and further use it to speed up the procedure of assembling fractional stiffness matrix.

Here we take equations with Riesz fractional derivatives as an example. Suppose we have a cubic
domain with a uniform grid in $\mathbb{R}^n$ and assume we have $N$ grid points on each
direction. As presented in \cref{table:sparsity}, when $n = 1$, the stiffness matrix is full.  If $n > 1$, the matrix is a sparse
matrix and the density of the matrix is proportional to $1/{N^{n-1}}$.  
\begin{table}[h]\centering\small
	\caption{Density of the stiffness matrix in different dimensions.}\label{table:sparsity}
	\begin{tabular}{cccc}\toprule
		Type             & 1-D      & 2-D          & 3-D          \\ \midrule
		Traditional Equ. & $O(1/N)$ & $O(1/{N^2})$ & $O(1/{N^3})$ \\
		Fractional Equ.  & $1$      & $O(1/{N})$   & $O(1/{N^2})$ \\ \bottomrule
	\end{tabular}
\end{table}

To accelerate the operation, we firstly construct a sparse matrix, called template matrix, of which the
structure is 'similar' with the result matrix $K$. 
{At the begining of assembling the fractional stiffness matrix, 
  we set $K = T$. And then we just add the element stiffness matrix 
  to $K$ by a C-routine which will save a lot of time. 
  After we finish the assembling process, we subtract $T$ from $K$. 
  Then we obtain the fractional stiffness matrix that we want.}
The readers can find the C routine at https://github.com/lrtfm/addsparse.

Now we briefly describe the idea of constructing the template matrix $T$.  Assume
$K_{\mathcal{N}\times\mathcal{N}}$ is a fractional stiffness matrix. Then we define template matrix
$T$ of $K$ as below
\begin{equation}
	T_{jk} = \left\{
	\begin{aligned}
		0, \quad \text{if}\ K_{jk} = 0,   \\
		1, \quad \text{if}\ K_{jk} \ne 0. 
	\end{aligned}
	\right.
\end{equation}
So, the process of constructing $T$ is actually to determine whether $K_{jk}$ is zero.
Here, we take  $\big( c_i(\mathbf{x}) \mathcal{D}^{\alpha}_i u(\mathbf{x}), {\mathcal{\hat D}^{\beta}_i} v(\mathbf{x}) \big)$ as an example,
i.e.~$K_{jk} = \big( c_i(\mathbf{x}) \mathcal{D}^{\alpha}_i \phi_k(\mathbf{x}), {\mathcal{\hat D}^{\beta}_i} \phi_j(\mathbf{x}) \big)$
where $ j,k \in \{1, 2,\cdots,\mathcal{N}\}$. As linear Lagrange element is used, we
assume the vertex corresponding to basis function $\phi_j(\mathbf{x})$ is $\mathbf{z}_j$ and define $\omega_{j}$
as the element patch of $\mathbf{z}_j$, i.e.~$\omega_{j} = \cup_{\mathbf{z}_j \in E} E$.
Define the lower and upper bound of $\omega_{j}$ in $\mathbf{e}_i$ direction as
\begin{equation}
	\left\{
	\begin{aligned}
		z_{j,min}^{(i)} = \min\{x_i| \mathbf{x} \in \omega_{j}\}, \\
		z_{j,max}^{(i)} = \max\{x_i| \mathbf{x} \in \omega_{j}\}. 
	\end{aligned}
	\right.
\end{equation}
According to the definition of fractional derivatives, if $[z_{j,min}^{(i)}, z_{j,max}^{(i)}] \cap [z_{k,min}^{(i)}, z_{k,max}^{(i)}] = \emptyset$, the element $K_{jk}$ and $K_{kj}$ must be zeros, i.e.~$T_{jk} = T_{kj} = 0$. Otherwise, $K_{jk}$ and $K_{kj}$ are non-zero, i.e~$T_{jk} = T_{kj} = 1$.

The efficiency of this algorithm can be seen from \cref{table:profile} which shows the time cost of MATLAB built-in function `plus' and C routine `addsparse' in assembling the
stiffness matrix of Riesz fractional derivatives. For fractional stiffness matrix, the acceleration effect of `addsparse' is remarkable. 

\begin{table}[!h]\small\centering
  \caption{Comparison of times cost by MATLAT {built-in} function `plus' and our C routine `addsparse'.}\label{table:profile}
	\begin{tabular}{crrcrr}\toprule
		Number of &    \multicolumn{2}{c}{$\alpha = 1$}   &  &  \multicolumn{2}{c}{$\alpha = 0.8$} \\ \cmidrule(lr){2-3}\cmidrule(lr){5-6}
		elements & `addsparse' & `plus' &   & `addsparse' & `plus'  \\
		\midrule
		940      & 0.06s       & 0.03s  &   & 0.07s       & 0.07s   \\
		2583     & 0.18s       & 0.13s  &   & 0.29s       & 0.62s   \\ 
		10746    & 1.05s       & 1.52s  &   & 1.73s       & 50.40s  \\    
		38139    & 7.52s       & 20.09s &   & 11.62s      & 945.53s \\
		\bottomrule
	\end{tabular}
\end{table}

\section{Application}\label{sec:examples}
In this section, we will use the methods proposed in previous section to solve the steady fractional diffusion 
equation with linear finite elements.

Let us consider FPDEs with variable coefficients in divergence form
\begin{equation}\label{equ:sfde2}
	\left\{
	\begin{aligned}
		  & \sum_{i=1}^{3} \mathcal{D}_i \Big(                     
		p_i(\mathbf{x}) \rfd{\beta_i}{a_i(\mathbf{x})}{x_i} u(\mathbf{x})
		- q_i(\mathbf{x}) \rfd{\beta_i}{x_i}{b_i(\mathbf{x})} u(\mathbf{x})
		\Big) = f(\mathbf{x}), 
		\quad \mathbf{x} \in \varOmega, \\
      & u(\mathbf{x}) = 0, \quad \mathbf{x} \in {\mathbb{R}^3\backslash\varOmega}, 
	\end{aligned}
	\right.
\end{equation}
where $p_i(x) >0, q_i(x) >0,\varOmega = \{\mathbf{x}| |\mathbf{x}| < r\}$ and
$r = 0.5$.
The variational form of this equation is $a(u, v) = (f, v)$, where
\begin{equation}
	\begin{aligned}
		a(u, v) & = - \sum_{i=1}^{3}\Big( 
		\big( p_i(\mathbf{x}) \rfd{\beta_i}{a_i(\mathbf{x})}{x_i} u, \mathcal{D}_i v \big)
		- \big( q_i(\mathbf{x}) \rfd{\beta_i}{x_i}{b_i(\mathbf{x})} u, \mathcal{D}_i v \big)
		\Big).
	\end{aligned}
\end{equation}

Denote  $\beta = (\beta_1, \beta_2, \beta_3)$ and set $(p_1, p_2, p_3) = (\cos(x_1), \cos(x_2), \cos(x_3))$, 
$(q_1, q_2, q_3) = (1-\cos(x_1), 1-\cos(x_2), 1-\cos(x_3))$. Let the exact solution be
\begin{equation}
	u(\mathbf{x}) = (x_1^2 + x_2^2 + x_3^2 - r^2)^2,
\end{equation} 
and $f(\mathbf{x})$ be the corresponding inhomogeneous term obtained by inserting $u(\mathbf{x})$ into the original equation.

\begin{figure}[!h]
	\centering
	\subfloat[Mesh with $h\approx 0.0894166$]{\includegraphics[width=0.34\textwidth]{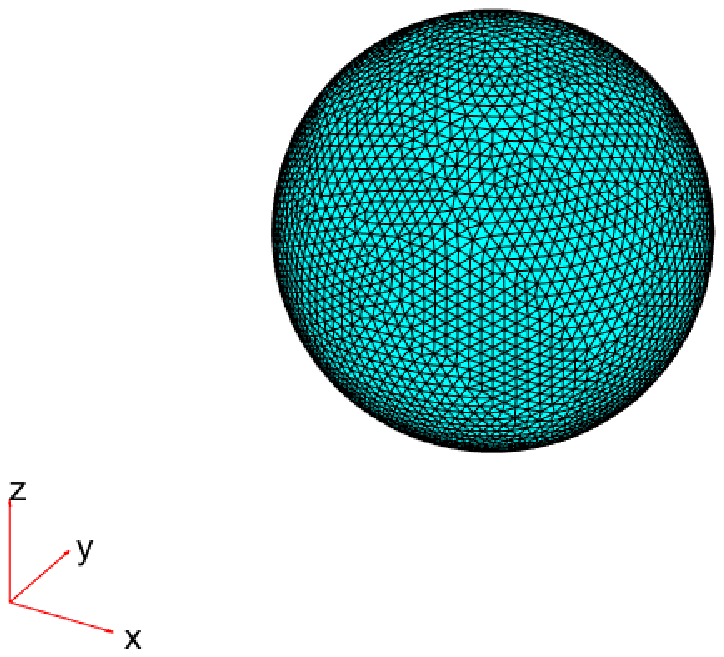}}
	\subfloat[Numerical result with $h\approx 0.0894166$]{\includegraphics[width=0.32\textwidth]{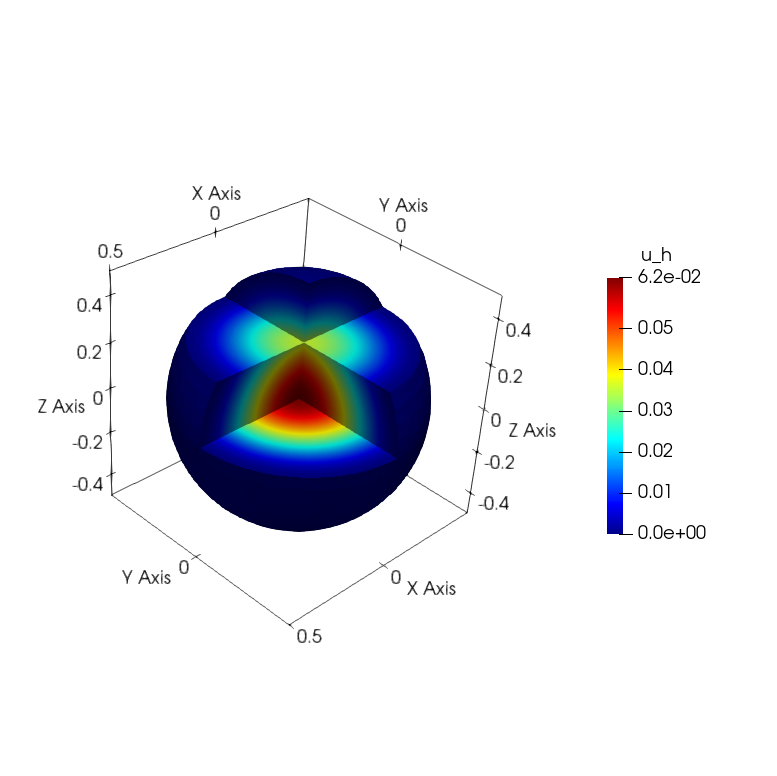}}
	\subfloat[Exact solution]{\includegraphics[width=0.32\textwidth]{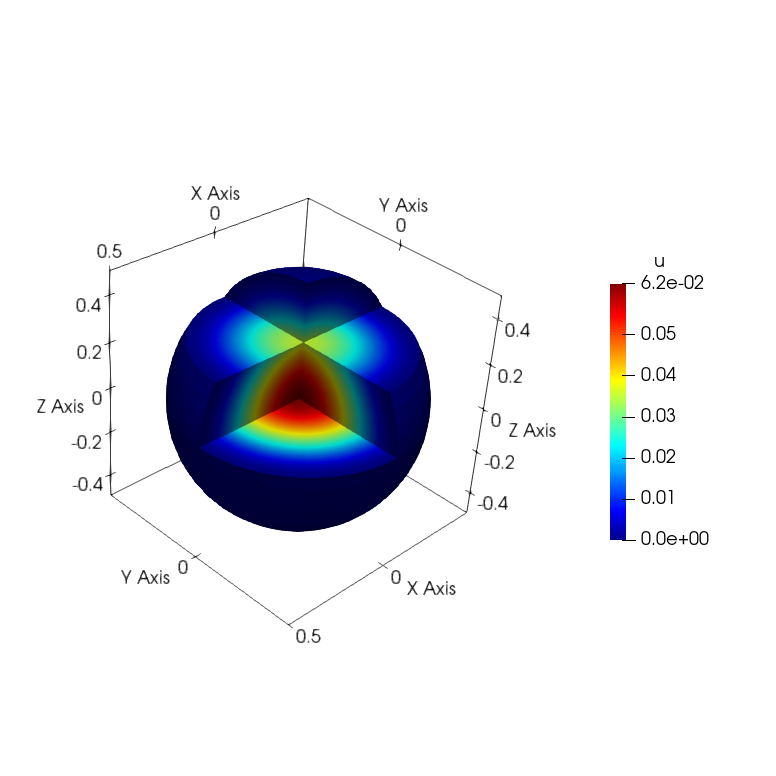}}
	\caption{Numerical result and exact solution for \eqref{equ:sfde2}.}
	\label{fig:gc2}
\end{figure}

The exact solution and the numerical result with $h \approx 0.0894166$ are presented in \cref{fig:gc2}.
The errors and convergence orders are shown in \cref{table:x1s} with $\beta = (0.8, 0.8, 0.8)$ and 
$\beta = (0.6, 0.7, 0.8)$, separately. According to \cite[Theorem 7.3]{ervinVariationalSolutionFractional2007}, the theoretical convergence order should be 2 which is consistent with our numerical results.
\begin{table}[!h]
	\centering\small
	\caption{$L^2$ and $L^\infty$ errors and convergence orders for \eqref{equ:sfde2}.}\label{table:x1s}
	\begin{tabular}{cccccccccc}
		\toprule
		\diagbox{$h$}{Errors}{$\beta$}  & \multicolumn{4}{c}{$(0.8, 0.8, 0.8)$} && \multicolumn{4}{c}{$(0.6, 0.7, 0.8)$} \\ \cline{2-5} \cline{7-10}
		& $\|u-u_h\|_{L^2}$ &   order & $\|u-u_h\|_{L^\infty}$ &  order  &
		& $\|u-u_h\|_{L^2}$ &   order & $\|u-u_h\|_{L^\infty}$ &  order \\
		\midrule
		0.270055  & 7.79e-04 &      & 2.24e-03 &      &   & 7.87e-04 &      & 2.30e-03 &      \\
		0.139416  & 2.64e-04 & 1.64 & 7.53e-04 & 1.65 &   & 2.74e-04 & 1.60 & 7.99e-04 & 1.60 \\
		0.0757345 & 8.05e-05 & 1.95 & 2.34e-04 & 1.91 &   & 8.68e-05 & 1.88 & 2.85e-04 & 1.69 \\
		\bottomrule
	\end{tabular}
\end{table}

\section{Conclusion}\label{sec:conclusion}
In this paper, we search the integration paths by a ray-simplex algorithm which can speed-up the procedures of assembling fractional stiffness matrices.  The algorithm is easy to implement and can be used in any $n$-dimensional space.  Furthermore, we apply the method to solved 3-D steady fractional problem. 
\begin{table}[h!]
	\centering
	\newcommand{\tabincell}[2][c]{\begin{tabular}{@{}#1@{}}#2\end{tabular}}
	\caption{Comparison of time used in solving fractional and integer elliptic equations}
	\label{table:comt}
	\begin{tabular}{rrccrc}
		\toprule
		& \multicolumn{2}{c}{Fractional elliptic equation} && \multicolumn{2}{c}{Elliptic equation}  \\
		\cmidrule(lr){2-3}\cmidrule(lr){5-6}
		\tabincell{Number \\of elements}    & \tabincell{Assemble \\ matrix}      & \tabincell{Solve \\ linear equations} &
		& \tabincell{Assemble \\ matrix}      & \tabincell{Solve \\ linear equations} \\ 
		\midrule
		4061   & 35.61s        & 0.02s &   & 0.92s  & 0.03s \\   
		64740  & 11m 34.60s    & 0.52s &   & 8.96s  & 0.09s \\
		268418 & 1h 36m 38.75s & 8.62s &   & 33.18s & 0.73s \\
		\bottomrule
	\end{tabular}
\end{table}

As shown in \Cref{table:comt}, although we propose a valid method on searching integral path and some {techniques} to speed up the addition of sparse matrix, the assembly of the fractional stiffness matrix is still much slower than that of the integer order stiffness matrix. It  is still needed to conduct in-depth research on fast algorithms of assembling fractional stiffness matrix on irregular domains.

\section*{Acknowledgements}

This research was supported by the National Natural Science Foundation of China
(Grant No.11601432 and No.11971386) and the Fundamental Research Funds for the 
Central Universities (Grant No.310201911cx025).

The authors would like to thank Professor Fawang Liu of Queensland University of Technology 
for his helpful comments on an earlier version of this manuscipt. The authors would like
to thank the anonymous reviewers for their useful comments, which have led to an 
improvement of the presentation.


\end{document}